\documentclass[12pt]{article}
\usepackage{geometry}
\usepackage{chngpage}
\usepackage{graphicx}
\usepackage{amsmath}
\usepackage{amsfonts}
\usepackage{amssymb}
\usepackage{latexsym}
\usepackage[brazil]{varioref}
\usepackage[english, brazil]{babel}
\usepackage[latin1]{inputenc}
\usepackage[dvips]{color}

\newtheorem{exemp}{Exemplo}[subsection]
\newtheorem{defn}{Defini\c{c}{\~a}o}[subsection]
\newtheorem{teo}{Teorema}[subsection]

\newtheorem{prop}[teo]{Proposição}
\newenvironment{dem}[1][Demonstra\c{c}{\~a}o]{\textbf{#1:}\ }
{\hfill\rule{1ex}{1ex}}

\begin{document}

\title{Espaços pseudo-topológicos e descrições do conceito de plausível}
\author{Tiago Augusto dos Santos Boza \thanks{ Email: boza.tiago@gmail.com. Pós em Filosofia, UNESP - FFC - Marília} \\
Hércules de Araujo Feitosa \thanks{ Email: haf@fc.unesp.br. Departamento de Matemática, UNESP - FC - Bauru}
 }

\maketitle \abstract{ A estrutura de espaço pseudo-topológico é uma variação do conceito de espaço topológico. Surgiu como uma formalização de um quantificador de primeira ordem não definível a partir dos usuais operadores existencial e universal, nomeado de o quantificador do plausível. Posteriormente, buscou-se uma formalização no contexto lógico proposicional, dada num sistema dedutivo e axiomático, com modelo algébrico. Apresentaremos estas concepções formais e como uma contribuição original, descrevemos a versão proposicional da lógica dos espaços quase topológicos em tableaux.

\medskip

\noindent{\bf Palavras Chave:} Espaços quase topológicos, Lógica modal, Sistema de axiomas, Modelos algébrico, Tableaux.

}


\section*{Introdução}

${}$\hspace{0,5cm} Grácio \cite{gra99} procurou formalizar quantificadores não lógicos, isto é, quantificadores distintos dos quantificadores $\exists$ e $\forall$, e que não pudessem ser definidos a partir destes usais quantificadores, mas que contemplassem aspectos do processo de generalização desenvolvido nos raciocínios indutivos. Estes quantificadores deveriam expressar proposições gerais, entretanto distintas do todo, do universal.

Para tanto, precisou de estruturas matemáticas que pudessem interpretar estes novos quantificadores com alguma sensatez intuitiva e interação com os demais aspectos da lógica.

Uma estrutura que surgiu nessa busca foi a de espaço pseudo-topológico, que é uma variação do usual conceito de espaço topológico. Esta estrutura foi idealizada como espaço de interpretação do quantificador ``para uma boa parte'' \cite{gra99}, o quantificador do plausível. Num primeiro momento, este conceito foi chamado de espaço topológico reduzido, mas posteriormente foi denominado de espaço  pseudo-topológico.

A partir desta estrutura, apresentamos a lógica do plausível de \cite{gra99}, que formaliza o conceito a partir de quantificadores estendidos, depois uma lógica proposicional com operador de caráter modal \cite{fng09}, motivada pela lógica do plausível, que formaliza aspectos do conceito de plausível no contexto proposicional e com modelos algébricos. Apesar da motivação, as duas versões de lógica têm flexões nas suas fundações.

Por outro lado, o método dos tableaux é baseado na refutação, de maneira que para se verificar a validade de uma fórmula $\psi$ em um sistema lógico, considera-se como hipótese a sua negação, ou seja, toma-se $\neg \psi$ e, então, utilizando
uma estrutura que se assemelha a uma árvore, são aplicadas as regras do sistema de tableaux. Baseamos nossas construções para tableaux nos textos \cite{smu71} e \cite{sfm06}.

Como uma contribuição original, apresentamos esta versão da lógica proposicional dos espaços quase topológicos numa versão de tableaux, a qual mostramos ser equivalente à versão axiomática de \cite{fng09}.


\section{Espaços pseudo-topológicos}

${}$\hspace{0,5cm} O conceito de espaço pseudo-topológico foi introduzido por Grácio \cite{gra99}, para interpretar o novo quantificador $P$, que gera sentenças do tipo $Px \varphi(x)$. Esta sentença tem o entendimento de que $\varphi(x)$ vale para uma boa parte dos indivíduos do domínio que interpretam a variável, ou ainda, conforme \cite{cag08}, que o conceito $\varphi(x)$ é ubíquo neste universo.

\begin{defn} Espaço pseudo-topológico é um par $(E, \Omega)$, em que $E$ é um conjunto não vazio e  $\Omega \subseteq \mathcal{P}(E)$ de maneira que:

$(E_1)$ se $A, B \in \Omega$, então $A \cap B \in \Omega$

$(E_2)$ se $A, B \in \Omega$, então $A \cup B \in \Omega$

$(E_3)$ $E \in \Omega$

$(E_4)$  $\emptyset \notin \Omega$. \end{defn}

\begin{defn} Os elementos do conjunto $\Omega$ são os abertos do espaço pseudo-topológico $(E, \Omega)$ e se $A^C \in \Omega$, então $A$ é um conjunto fechado no espaço $(E, \Omega)$. \end{defn}

\begin{exemp} Sejam $E \neq \emptyset$ e $\Omega = \{E\}$.

Certamente  $\emptyset \notin \Omega$, mas $E \in \Omega$. As condições (i) e (ii) são trivialmente satisfeitas. Então, $(E, \Omega)$ é um espaço pseudo-topológico. \end{exemp}

\begin{exemp} Seja $E \neq \emptyset$. Para $a \in E$, seja $\Omega = \{B \subseteq E : a \in B \}$. 

Então, de modo óbvio verificamos que $\emptyset \notin \Omega$, mas $E \in \Omega$. 
E dessa forma obtemos que $(E, \Omega)$ é um espaço pseudo-topológico. \end{exemp}

\begin{prop} \label{cadeia} Se $(E, \Omega)$ é um espaço pseudo-topológico, então dois quaisquer abertos de $\Omega$ não são disjuntos. \\
\begin{dem} Se $A$ e $B$ são dois abertos disjuntos, pelo axioma (i),  $A \cap B \in \Omega$, mas isto contradiz o axioma (iv), pois neste caso $A \cap B = \emptyset$. \end{dem}
\end{prop}

Assim, conjuntos disjuntos não podem ser simultaneamente abertos. Segue daí, que não podem existir dois conjuntos unitários distintos e abertos.


\subsection{A motivação quantificacional}

${}$\hspace{0,5cm} A formalização lógica quantificacional do espaços pseudo-topológicos surgiu em \cite{gra99} como uma extensão da lógica clássica de primeira ordem com igualdade $\mathcal{L}$. Para detalhes sobre $\mathcal{L}$ ver \cite{eft84} ou \cite{fep05}. Foi denominada por Grácio de lógica do plausível e denotada por $\mathcal{L}(P)$. \\

Dada $\mathcal{L}$, a lógica estendida $\mathcal{L}(P)$ é determinada pelos seguintes acréscimos: \\

- a linguagem $L$ de $\mathcal{L}(P)$ conta com um novo símbolo de quantificador $P$ e sentenças do tipo $Px \varphi(x)$ são bem formadas em  $\mathcal{L}(P)$. \\

- axiomas específicos do quantificador $P$:

$(A_1)$ $Px \varphi(x) \land Px \psi(x) \to Px(\varphi(x) \land Px \psi(x))$

$(A_2)$ $Px \varphi(x) \land Px \psi(x) \to Px(\varphi(x) \lor Px \psi(x))$

$(A_3)$  $\forall x \varphi(x) \rightarrow Px \varphi(x)$

$(A_4)$ $Px \varphi(x) \rightarrow \exists x \varphi(x)$

$(A_5)$ $\forall x (\varphi(x) \to \psi(x)) \to (Px \varphi(x) \to Px \psi(x))$

$(A_6)$ $Px \varphi(x) \to Py \varphi(y)$, se $y$ é livre para $x$ em $\varphi(x)$. \\

- regradas de dedução:

(MP) \emph{Modus Ponens}: $\varphi, \varphi \to \psi \vdash \psi$

(Gen) Generalização: $\varphi \vdash \forall x \varphi(x)$. \\

Os demais conceitos sintáticos usuais como sentença, demonstração, teorema, dedução, consistência e outros são definidos do modo padrão. \\

As estruturas adequadas para $\mathcal{L}(P)$, denominadas por Grácio de estrutura do plausível, também são extensões das estruturas de primeira ordem $\mathcal{A}$. \\

Assim dada uma estrutura $\mathcal{A}$, consideremos que o seu domínio é denotado por $A$. Uma estrutura do plausível, denotada por $\mathcal{A}^{\Omega}$, é determinada a partir $\mathcal{A}$ pelo acréscimo de um espaço pseudo-topológico $\Omega$ sobre o universo $A$.

A interpretação dos símbolos de relação, função e constante é a mesma de $\mathcal{L}$ com relação à $\mathcal{A}$.

\begin{defn} A satisfação de uma sentença do tipo $Px \varphi(x)$ em $\mathcal{A}^{\Omega}$ e definida indutivamente por:

- se $\varphi$ é uma fórmula cujas variáveis livres estão em $\{x\} \cup \{y_1, ..., y_n\}$ e $\overline{a} =
(a_1, ..., a_n)$ é uma sequência de elementos de $A$, então: 
$$\mathcal{A}^{\Omega} \vDash Px \varphi[x, \overline{a}] \Leftrightarrow \{b \in A : \mathcal{A}^{\Omega} \vDash
[b, \overline{a}]\} \in \Omega.$$
\end{defn} 

Da maneira usual, para a sentença $Px \varphi(x)$:
$$\mathcal{A}^{\Omega} \vDash Px \varphi(x) \Leftrightarrow \{a \in A : \mathcal{A}^{\Omega} \vDash \varphi(a)\} \in \Omega.$$

As outras noções semânticas como modelo, validade, implicação lógica, entre outras, são apropriadamente adaptadas a partir da interpretação de $\mathcal{L}$ em $\mathcal{A}$. \\

Grácio \cite{gra99} provou que as estruturas do plausível são modelos corretos e completos para $\mathcal{L}(P)$.


\subsection{Uma formalização axiomática e proposicional}

${}$\hspace{0,5cm} Apresentamos a versão proposicional da lógica do plausível, conforme \cite{fng09}, que procura formalizar os aspectos de um espaço pseudo-topológico no contexto lógico proposicional, sem a presença dos quantificadores, inclusive o do plausível, que é substituído por um operador unário. \\

Esta lógica será denotada por $\mathbb{L}(\nabla)$. Como mencionado, ela estende a lógica proposicional clássica (LPC) na linguagem $L(\neg, \land, \lor, \to)$ com o acréscimo do operador unário $\nabla$, donde obtermos a linguagem proposicional $L(\neg, \land, \lor, \to, \nabla)$.

A lógica fica determinada pelo seguinte: \\

- Axiomas:

$(Ax_0)$ LPC

$(Ax_1)$ $(\nabla \varphi \land \nabla \psi) \to \nabla(\varphi \land \psi)$

$(Ax_2)$ $\nabla(\varphi \lor \neg \varphi)$ 

$(Ax_3)$ $\nabla \varphi \to \varphi$. \\

- Regras de dedução:

$(MP)$ \emph{Modus Ponens}

$(R\nabla)$ $\vdash \varphi \to \psi\ / \vdash \nabla \varphi \to \nabla \psi$. \\

A intuição para este operador de plausível é de algo que pode ser explicado numa teoria sem provocar inconsistência. Assim, o plausível não se assemelha ao possível, pois podemos ter: ``é possível que chova amanhã'' e ``é possível que não chova amanhã'', porém ``não é possível que chova e não chova amanhã''.

O conceito está vinculado, por exemplo, com a existência de uma demonstração de um fato numa teoria consistente. Assim, não pode haver uma prova de $\varphi$ e uma outra de $\neg \varphi$.

Diante disso, $\varphi$ e $\psi$ são plausíveis se, e somente se, $\varphi \land \psi$ é plausível. Toda tautologia é plausível e se $\varphi$ é plausível, então vale a proposição $\varphi$. A regra $(R\nabla)$ diz que se há uma prova de $\varphi \to \psi$ e $\varphi$ é plausível, então também $\psi$ é plausível. \\

	Pode ser demonstrado o seguinte.
	
\begin{prop} (i) $\vdash \neg \nabla \bot$

(ii) $\vdash \nabla \varphi \to \nabla (\varphi \lor \psi)$

(iii) $\vdash \varphi \Rightarrow\ \vdash \nabla \varphi$

(iv) $\vdash \varphi \to \neg \nabla \neg \varphi$

(v) $\vdash \nabla \varphi \to \neg \nabla \neg \varphi$

(vi) $\vdash \nabla \neg \varphi \to \neg \nabla \varphi$.
\end{prop}

\begin{prop} \label{ou} $\nabla \varphi \to \nabla(\varphi \lor \psi) \Leftrightarrow (\nabla \varphi \lor \nabla \psi) \to \nabla(\varphi \lor \psi)$. \\
\begin{dem} $(\Rightarrow)$ Da hipótese, $\nabla \varphi \to \nabla(\varphi \lor \psi)$ e $\nabla \psi \to \nabla(\varphi \lor \psi)$. Daí, $(\nabla \varphi \lor \nabla \psi) \to \nabla(\varphi \lor \psi)$.

$(\Leftarrow)$ Como $\nabla \varphi \to (\nabla \varphi \lor \nabla \psi)$, então segue da hipótese que $\nabla \varphi \to \nabla(\varphi \lor \psi)$. \end{dem}
\end{prop}


\subsection{Álgebras dos espaços pseudo-topológicos}

${}$\hspace{0,5cm} Agora, apresentamos a álgebra do plausível, que correspondem à versão algébrica da lógica da seção anterior, conforme \cite{fng09}.

\begin{defn} Álgebra do Plausível é uma estrutura $\mathbb{P} = (P, 0, 1, \land, \lor, \sim, \sharp)$, em que $(P, 
0, 1, \land, \lor, \sim)$ é uma álgebra de Boole e $\sharp$ é o operador do plausível, sujeito a: 

$(a_1)$  $\sharp a \land \sharp b \leq \sharp (a \land b)$

$(a_2)$  $\sharp a \leq \sharp (a \lor b)$

$(a_3)$  $\sharp a \leq a$

$(a_4)$  $\sharp 1 = 1$. 	\end{defn}

\begin{defn} Um elemento $a \in P$ é plausível $a \neq 0$ e $\sharp a = a$.	\end{defn}  

Embora $\sharp 0 = 0$, por definição, $0$ não é plausível.

\begin{prop} Se $\mathbb{P} = (P, 0, 1, \land, \lor, \sim, \sharp)$ é uma álgebra do plausível e $a, b \in P$, então:

(i) $\sharp a \leq \sharp(a \lor b)$

(ii) $a \leq b \Rightarrow \sharp a \leq \sharp b$

(iii) $\sharp a \lor \sharp b \leq \sharp (a \lor b)$. 
\end{prop}

\begin{prop} Para cada álgebra do plausível $\mathbb{P} = (P, 0, 1, \land, \lor, \sim, \sharp)$ existe um monomorfismo $h$  de $P$ num espaço pseudo-topológico de conjuntos definidos em $\mathcal{P}(\mathcal{P}(P))$.
\end{prop}

Em \cite{fng09} há uma demonstração da adequação de $\mathbb{L}(\nabla)$ com relação às álgebras do plausível $\mathbb{P}$. \\

A interpretação de $\mathbb{L}(\nabla)$ numa pseudo-topologia $(E, \Omega)$ é uma função $v$ com as seguintes características:
	
As sentenças universais que valem para todos os indivíduos, são agora tomadas pelas tautologias, que são interpretados no universo $E$ da pseudo-topologia. 

Se $\varphi$ e $\psi$ estão $\Omega$ e, portanto, são ubíquos, então o mesmo vale  $\varphi \land \psi$.

Se $\varphi$ ou $\psi$ está $\Omega$, então o mesmo vale  $\varphi \lor \psi$. Esta condição implica que se $\varphi, \psi \in \Omega$, então $\varphi \lor \psi \in \Omega$, mas é um pouco mais básica. Aqui está uma variação entre os dois sistemas.

Se $\varphi$ e $\psi$ são equivalentes, então uma delas está numa pseudo-topologia se, e somente se, também a outra está.

A interpretação de $\bot$ não pode estar $\Omega$. Para tanto, usamos o axioma $(Ax_3)$ $\nabla \bot \to \bot$ que implica esta condição no contexto proposicional. Mas também tem mais exigência que $(E_4)$. \\

Talvez possamos, em algum momento, refinar estas variações, mas por ora assumamos a lógica $\mathbb{L}(\nabla)$.


\section{Tableaux para a lógica $\mathbb{L}(\nabla)$}

${}$\hspace{0,5cm} Tendo em conta os sistemas de tableaux de \cite{smu71} e \cite{sfm06} introduziremos um sistema de tableaux para $\mathbb{L}(\nabla)$. Daremos as regras dos tableaux e para mostrarmos a equivalência entre os dois sistemas dedutivos, mostraremos que se $\varphi$ é um teorema de  $\mathbb{L}(\nabla)$, então o seu tableau fecha e, por outro lado, se o tableau da fórmula $\varphi$ fecha, então ela é um teorema de $\mathbb{L}(\nabla)$.

O sistema de tableaux a ser desenvolvido nesta seção deve ser algorítmico e, tanto quanto possível, preservar as características dos tableaux clássicos. Deve derivar exatamente o mesmo que o sistema axiomático $\mathbb{L}(\nabla)$, nem mais, nem menos. 

Denotaremos por $\mathcal{T}_{L(P)}$ a lógica proposicional do plausível na versão de tableaux.

\begin{defn} A linguagem do sistema $\mathcal{T}_{L(P)}$ é determinada a partir da
formalização dos tableaux originais apresentados por \cite{smu71}, acrescida apenas dos seguintes itens:\\
(i) O alfabeto de $\mathcal{T}_{L(P)}$, denotado por $(Alf_(T_Pl))$, é constituído pelos itens
apresentados por \cite{smu71}, acrescido do operador $\nabla$;\\
(ii) O conjunto de fórmulas de $\mathcal{T}_{L(P)}$, denotado por $(For_(T_Pl))$, é dado recursivamente pelas fórmulas apresentadas originalmente, acrescido da seguinte cláusula:

- Se $A$ é uma fórmula, então $\nabla A$ também é uma fórmula de $\mathcal{T}_{L(P)}$; \\
(iii) O conjunto de regras de dedução do sistema $\mathcal{T}_{L(P)}$ é formado pelas
regras originais de dedução dos tableaux para a LPC, acrescido das regras de dedução \textit{específicas} para o operador $\nabla$, que serão introduzidas na definição posterior; \\
(iv) Um ramo do tableau fecha ao obtermos $A$ e $\neg A$ no mesmo ramo ou ainda o símbolo de contradição, $\bot$.

\end{defn}

\begin{defn} As regras de expansão de $\mathcal{T}_{L(P)}$ para o operador $\nabla$ são as
seguintes: \\
$(R_1)$ $\dfrac{\nabla A}{A}$	\\
A motivação para esta regra vem, naturalmente, do axioma $(Ax_3)$. \\
\\$(R_2)$ $\dfrac{(\neg \nabla A) \land (\Vdash A)}{\bot}$ \\
\\Esta regra nos permite a obtenção da validade do axioma $(Ax_4)$. Bem como, esta regra nos permite o caráter algorítmico que pretendemos em nosso sistema. 
Pois, quando encontramos uma expressão $\neg \nabla A$, então testamos $A$. Se a fórmula $A$ é válida, isto é, o tableau de $\neg A$ fecha, então incluímos $\bot$ e fechamos o tableau, pois se $A$ é válida, então tem que valer $\nabla A$. Agora, se $A$ não é válida, então expandimos o
tableau com uma das regras seguintes: \\
\\$(R_3)$ $\dfrac{\neg \nabla (A \land B)}{\neg \nabla A \lor \neg \nabla B}$ \\
\\$(R_4)$ $\dfrac{\neg \nabla (A \lor B)} {\neg \nabla A \land \neg \nabla B}$ \\
\\As regras $(R_3)$ e $(R_4)$ dão conta dos axiomas $(Ax_1)$ e $(Ax_2)$. \\
$(R_{5A})$ $\dfrac{\neg (A \rightarrow B)}{\neg \nabla (\neg A \lor B)}$ \\
$(R_{5B})$ $\dfrac{\neg (A \leftrightarrow B)}{\neg \nabla (A \rightarrow B) \land (B \rightarrow A)}$ \\
\\As regras $(R_{5A})$ e $(R_{5B})$ tem apenas a função de dar caráter algorítmico ao tableau ao dizer exatamente qual caminho seguir.\\
\\$(R_6)$ $\dfrac{(A \leftrightarrow B)} {(\nabla A \land \nabla B) \lor (\neg \nabla A \land \neg \nabla B)}$ \\
\\Finalmente, a regra $(R_6)$ tem a incumbência de validar a regra $(R \nabla)$
que apenas se aplica sobre bicondicionais sabidamente válidas.

\end{defn}

Agora, mostraremos alguns resultados importantes, por tableaux, demonstrados na versão axiomática da lógica proposicional do plausível.

\begin{prop} Nenhuma contradição é plausível.\\
\begin{dem}
Para mostrarmos o que queremos basta demonstrarmos que o tableau para $\neg \nabla (A \land \neg A)$ fecha.
\begin{center}
 $\neg \neg \nabla (A \land \neg A)$ \\
 $\nabla (A \land \neg A)$ \\
 $(A \land \neg A)$ \\
 $A$ \\
 $\neg A$ \\
 $_x$
 \end{center} 
Portanto, demonstramos o que queríamos.
\end{dem}
\end{prop}

\begin{prop} $\Vdash A \Rightarrow \Vdash \nabla A$ \\
\begin{dem}
\begin{center}
$\Vdash A$ \\
$\neg \nabla A$ \\
$\bot$ \\
$_x$
\end{center}
Portanto,$\Vdash A \Rightarrow \Vdash \nabla A$.
\end{dem}
\end{prop}

\begin{prop} $\Vdash \nabla(\nabla A \Rightarrow \nabla A)$ \\
\begin{dem}
\begin{center}
$\neg \nabla(\nabla A \Rightarrow \nabla A)$ \\
$\nabla(\nabla A \Rightarrow \nabla A)$ \\
$\bot$ \\
$_x$
\end{center}
Portanto,$\Vdash \nabla(\nabla A \Rightarrow \nabla A)$\end{dem}
\end{prop}

A seguir, mostraremos a equivalência do sistema Hilbertiano para a lógica proposicional do plausível, com relação ao sistema de tableaux apresentado neste capítulo. Mais especificamente, os resultados de correção e completude.

\begin{prop}
 $\textit{(Correção)}$ $\Gamma \vdash A \rightarrow \Gamma \Vdash A$. \\
 \begin{dem} Façamos essa demonstração por indução sobre o comprimento da dedução. Para n = 1, temos que $A$ pertence a $\Gamma$ ou $A$ é um axioma:
(i) Se $A \in \Gamma$ ou $A$ é um dos axiomas da LPC, nada temos a demonstrar,
pois foi demonstrado anteriormente, por $\cite{smu71}$.
(ii) Seja $A$ é um dos axiomas da lógica proposicional do plausível:
Mostraremos, relativo a cada axioma, que existe um tableau fechado
para $\Gamma \cup \neg A$: \\
\\- Para $(Ax_1)$, temos $(\nabla A \land \nabla B) \rightarrow \nabla (A \land B)$:
\begin{center}
$\Gamma$ \\
$\nabla A \land \nabla B$ \\
$\neg \nabla(A \land B)$ \\
$\nabla A$ \\
$\nabla B$ \\
$\curlywedge$ \\
$\neg \nabla A$ ${}$\hspace{1,0 cm} $\neg \nabla B$ \\
$_x$ ${}$\hspace{1,3cm} $_x$ \\
\end{center} 
- Para $(Ax_2)$, temos: $(\nabla A \lor \nabla B) \rightarrow \nabla (A \lor B)$:
\begin{center}
$\Gamma$ \\
$\nabla A \lor \nabla B$ \\
$\neg \nabla (A \lor B)$ \\
$\neg \nabla A$ \\
$\neg \nabla B$ \\
$\curlywedge$ \\
$\nabla A$ ${}$\hspace{1,0 cm} $\nabla B$ \\
$_x$ ${}$\hspace{1,3 cm} $_x$ \\
\end{center} 
- Para $(Ax_3)$, temos: $\nabla A \rightarrow A$:
\begin{center}
$\Gamma$ \\
$\nabla A$ \\
$\neg A$ \\
$A$ \\
$_x$\\
\end{center} 
- Para $(Ax_4)$, temos: $\nabla (A \lor \neg A)$:
\begin{center}
$\Gamma$ \\
$\neg \nabla (A \lor \neg A)$ \\
$\Vdash (A \lor \neg A)$ \\
$\bot$ \\
$_x$ \\
\end{center} 
Agora, tomemos como Hipótese de Indução que: $\Gamma \vdash A_n \Rightarrow \Vdash A_n$, para todo $n \leqslant k$. Há três possibilidades para o passo seguinte, $k+1$, da indução: \\

i) $A_{k+1}$ é uma premissa;

ii) $A_{k+1}$ é um esquema de axioma da lógica proposicional do
plausível;

iii) $A_{k+1}$ é deduzida a partir da regra $(R \nabla)$, ou seja, $\vdash A \leftrightarrow B \Rightarrow \vdash \nabla A \leftrightarrow \nabla B$. \\
\\Para os itens (i) e (ii), nada temos a demonstrar, pois essa demonstração é a mesma feita no caso base da indução. Já, para o caso (iii), é preciso que analisemos a regra $(R \nabla)$. Suponhamos que há uma dedução de $A_k \equiv A \leftrightarrow B$ a partir de $\Gamma$ e que $A_{k+1} \equiv \nabla A \leftrightarrow \nabla B$.
Segue daí que $\Gamma \vdash A \leftrightarrow B$ e, pela hipótese de indução, que $\Gamma \Vdash A \leftrightarrow B$.
Desse modo, basta mostrarmos que o tableau para $\Gamma \cup ${$ \neg (\nabla A \leftrightarrow \nabla B) $}$ $ também fecha, ou seja, que $\Gamma \Vdash \nabla A \leftrightarrow \nabla B$.
Agora, construindo o tableau que queremos, temos:
\begin{center}
$\Gamma$ \\
$\neg (\nabla A \leftrightarrow \nabla B)$ \\
$\curlywedge$ \\
$\neg (\nabla A \rightarrow \nabla B)$ ${}$\hspace{1,5 cm} $\neg (\nabla B \rightarrow \nabla A)$ \\
$\nabla A$ ${}$\hspace{3,0 cm} $\nabla B$ \\
$\neg \nabla B$ ${}$\hspace{3,0 cm} $\neg \nabla A$ \\
$\curlywedge$ ${}$\hspace{3,0 cm} $\curlywedge$ \\
$\nabla A$ ${}$\hspace{0,75 cm}	$\neg \nabla A$ ${}$\hspace{0,75 cm} $\nabla A$ ${}$\hspace{0,75 cm} $\neg \nabla A$ \\
$\nabla B$ ${}$\hspace{0,75 cm}	$\neg \nabla B$ ${}$\hspace{0,75 cm} $\nabla B$ ${}$\hspace{0,75 cm} $\neg \nabla B$ \\
$_x$ ${}$\hspace{1,4 cm}	$_x$ ${}$\hspace{1,4 cm} $_x$ ${}$\hspace{1,4 cm} $_x$ \\
\end{center}
Usamos a hipótese de indução para a última parte do tableau.
Com esses resultados concluímos que o sistema $\mathcal{T}_{L(P)}$ é correto, isto é, $\Gamma \vdash A \Rightarrow \Gamma \Vdash A$.
\end{dem}
\end{prop}

Para demonstrarmos a recíproca desta dedução no sistema $\mathcal{T}_{L(P)}$ com
relação ao sistema hilbertiano utilizaremo-nos das regras denominadas de $\textit{regras de inferência}$. Pois, segundo $\cite{fep05}$, essas regras são argumentos válidos, ou seja, conduzem premissas verdadeiras em conclusões verdadeiras. \\

\begin{prop} $\textit{(Completude)}$ $\Gamma \Vdash A \Rightarrow \Gamma \vdash A$. \\
\begin{dem} 
Assim como foi feito na demonstração anterior, estenderemos a demonstração feita para a lógica proposicional clássica para o operador $\nabla$. Mais especificamente, para cada regra de $\mathcal{T}_{L(P)}$ será obtida uma dedução, muitas delas usando o princípio da redução ao absurdo $(RAA)$, isto é, uma dedução indireta na lógica proposicional do plausível. \\
\\ - Para a regra $(R_1)$, temos: \\
1. $\nabla A$ ${}$\hspace{5,3 cm}   p.\\
2. $\neg A$ ${}$\hspace{5,3 cm} p.p. \\
3. $\nabla A \rightarrow A$ ${}$\hspace{4,3 cm} $(Ax_3)$ \\
4. $A$ ${}$\hspace{5,5 cm} Modus Ponens em 1 e 3 \\
5. $A \land \neg A$ ${}$\hspace{4,5 cm} Conjunção em 2 e 4 \\
6. $A$ ${}$\hspace{5,5 cm} RAA de 2 a 5 \\
\\ - Para a regra $(R_2)$, temos: \\
1. $\neg \nabla A$ ${}$\hspace{5,1 cm} p. \\
2. $\vdash A$ ${}$\hspace{5,3 cm} p. \\
3. $\vdash A \leftrightarrow (B \lor \neg B)$ ${}$\hspace{3,0 cm} Equivalência de Teoremas \\
4. $\vdash \nabla A \leftrightarrow \nabla (B \lor \neg B)$ ${}$\hspace{2,4 cm} $(R \nabla)$ em 3 \\
5. $\vdash \nabla (B \lor \neg B)$ ${}$\hspace{3,6 cm} $(Ax_4)$ \\
6. $\vdash \nabla (B \lor \neg B) \rightarrow \nabla A$ ${}$\hspace{2,4 cm} CPC em 4 \\
7. $\vdash \nabla A$ ${}$\hspace{4,9 cm} Modus Ponens em 5 e 6 \\
8. $\nabla A \land \neg \nabla A$ ${}$\hspace{3,9 cm} Conjunção em 1 e 7 \\
9. $\bot$ ${}$\hspace{5,5 cm} Equivalência em 8 \\
\\Esta regra $(R_2)$ é um pouco estranha. As regras devem ter a característica de conduzirem sentenças válidas em sentenças válidas, mas neste caso chegamos a uma contradição. Isto ocorre porque as premissas não podem concomitantemente ocorrerem. As regras devem ser tais que se as premissas
são válidas então a conclusão é válida. Neste caso, as premissas não são
ambas válidas, logo a conclusão é contraditória. \\
\\ - Para a regra $(R_3)$, temos: \\
1. $\neg \nabla (A \land B)$ ${}$\hspace{5,0 cm} p. \\
2. $\neg (\neg \nabla A \lor \neg \nabla B)$ ${}$\hspace{4,1 cm} p.p. \\
3. $\nabla A \land \nabla B$ ${}$\hspace{5,1 cm} De Morgan em 2 \\
4. $\nabla A \land \nabla B \rightarrow \nabla (A \land B)$ ${}$\hspace{2,9 cm} $(Ax_1)$ \\
5. $\nabla (A \land B)$ ${}$\hspace{5,1 cm} Modus Ponens em 3 e 4 \\
6. $(\nabla(A \land B)) \land (\neg \nabla (A \land B))$ ${}$\hspace{2,15 cm} Conjunção em 1 e 5 \\
7. $\neg \nabla (A \land B) / \neg \nabla A \lor \neg \nabla B$ ${}$\hspace{2,4 cm} RAA de 1 a 6 \\
\\ - Para a regra $(R_4)$, temos: \\
1. $\neg \nabla (A \lor B)$ ${}$\hspace{5,0 cm} p. \\
2. $\neg (\neg \nabla A \land \neg \nabla B)$ ${}$\hspace{4,15 cm} p.p. \\
3. $\nabla A \lor \nabla B$ ${}$\hspace{5,15 cm} De Morgan em 2 \\
4. $(\nabla A \lor \nabla B) \rightarrow \nabla (A \lor B)$ ${}$\hspace{2,6 cm} $(Ax_2)$ \\
5. $\nabla (A \lor B)$ ${}$\hspace{5,15 cm} Modus Ponens em 3 e 4 \\
6. $(\neg \nabla (A \lor B)) \land (\nabla (A \lor B))$ ${}$\hspace{2,25 cm} Conjunção em 1 e 5 \\
7. $\neg \nabla (A \lor B) / \neg \nabla A \land \neg \nabla B$ ${}$\hspace{2,5 cm} RAA de 1 a 6 \\
\\ - Para a regra $(R_{5A})$, temos: \\
1. $\neg \nabla (A \rightarrow B)$ ${}$\hspace{4,8 cm} p. \\
2. $\neg \neg \nabla (\neg A \lor B)$ ${}$\hspace{4,4 cm} p.p. \\
3. $\nabla (\neg A \lor B)$ ${}$\hspace{4,85 cm} Dupla Negação em 2 \\
4. $\nabla (A \rightarrow B)$ ${}$\hspace{4,9 cm} Equivalência em 3 \\   
5. $(\neg \nabla (A \rightarrow B)) \land (\nabla (A \rightarrow B))$ ${}$\hspace{1,8 cm} Conjunção em 1 e 4 \\
6. $\neg \nabla (A \rightarrow B) / \neg \nabla (\neg A \lor B)$ ${}$\hspace{2,3 cm} RAA de 1 a 5 \\
\\ - Para a regra $(R_{5B})$, temos: \\
1. $\neg \nabla (A \leftrightarrow B)$ ${}$\hspace{6,2 cm} p. \\
2. $\neg \neg \nabla ((A \rightarrow B) \land (B \rightarrow A))$ ${}$\hspace{3,7 cm} p.p. \\
3. $\nabla ((A \rightarrow B) \land (B \rightarrow A))$ ${}$\hspace{4,2 cm} Dupla Negação em 2 \\
4. $\nabla (A \leftrightarrow B)$ ${}$\hspace{6,35 cm} Equivalência em 3 \\
5. $(\neg \nabla (A \leftrightarrow B)) \land (\nabla (A \leftrightarrow B))$ ${}$\hspace{3,25 cm} Conjunção em 1 e 4 \\
6. $\neg \nabla (A \leftrightarrow B) / \neg \nabla ((A \rightarrow B) \land (B \rightarrow A))$ ${}$\hspace{1,55 cm} RAA de 1 a 5 \\
\\ - Para a regra $(R_6)$, temos: \\
1. $\vdash A \leftrightarrow B$ ${}$\hspace{6,7 cm} p. \\
2. $\vdash \neg (\nabla A \land \nabla B)$ ${}$\hspace{5,60 cm} p.p. \\
3. $\vdash \neg (\neg \nabla A \land \neg \nabla B)$ ${}$\hspace{5,10 cm} p.p. \\
4. $\vdash \neg (\nabla A \leftrightarrow \nabla B)$ ${}$\hspace{5,40 cm} CPC em 2 e 3 \\
5. $\vdash \nabla A \leftrightarrow \nabla B$ ${}$\hspace{5,95 cm} $(R \nabla)$ em 1 \\
6. $\vdash (\nabla A \leftrightarrow \nabla B) \land \neg(\nabla A \leftrightarrow \nabla B)$ ${}$\hspace{2,8 cm} Conjunção em 5 e 6 \\
7. $\vdash (\nabla A \land \nabla B) \lor (\neg \nabla A \land \neg \nabla B)$ ${}$\hspace{2,9 cm} RAA de 2 a 6 \\

Com esses resultados concluímos que $\Gamma \Vdash A \Rightarrow \Gamma \vdash A$.
\end{dem}
\end{prop}

\begin{prop} $\textit{Adequação Forte}$ $\Gamma \Vdash A \Leftrightarrow \Gamma \vdash A \Leftrightarrow \Gamma \models A$. \\
\begin{dem}
Segue imediatamente das proposições anteriores e Corollary
5.5 de $\cite{fng09}$. 
\end{dem}
\end{prop}

Da proposição acima concluímos que o presente artigo estabelece a equivalência entre a lógica proposicional do plausível e o sistema de tableaux $\mathcal{T}_{L(P)}$. 


\section*{Considerações finais}

${}$\hspace{0,5cm} O presente artigo mostra a equivalência entre a lógica proposicional
do plausível e o sistema de tableaux denotado por $\mathcal{T}_{L(P)}$. \\
O sistema $\mathcal{T}_{L(P)}$ procurou espelhar a lógica proposicional do plausível
através de seus axiomas, que neste caso, resgatam muito dos seus modelos.
Contudo, não temos claro que seja o único, nem o melhor sistema de
tableaux para esta lógica.

A tradição de tableaux para lógicas modais usa com bastante
frequência variações nas relações dadas nos seus modelos, pois também,
com frequência tratam de lógicas modais normais, que admitem o axioma K
e, daí, os modelos de Kripke. Porém a lógica proposcional do plausível seria,
no contexto modal, uma lógica subnormal e os modelos de Kripke não se
aplicam de modo imediato. Assim, julgamos conveniente olharmos para os
axiomas e buscar regras que os replicassem em $\mathcal{T}_{L(P)}$.

Agora, Em geral, um sistema de tableaux é considerado eficiente por ser
baseado no princípio das subfórmulas. No entanto, nesse trabalho este
princípio não parece estar integralmente presente. Mesmo assim, tentamos
de uma certa forma obedecermos tal princípio, por exemplo, ao tentarmos
uma regra que naturalmente viria do axioma $(Ax_3)$ pensaríamos, pela contra-
positiva, na seguinte regra: $\neg \nabla A / \neg A$, mas não nos pareceu muito intuitivo
do ponto de vista deste princípio. Desse modo, foi dada uma regra,
naturalmente equivalente, à saber, a regra $(R_1)$.

Diante disso, podemos dizer que um dos intuitos deste artigo seria a obtenção de
um sistema de dedução para a lógica proposicional do plausível
algoritmicamente mais eficiente do que a abordagem hilbertiana na qual ela
foi inicialmente apresentada. E isto é possível através do sistema de
tableaux $\mathcal{T}_{L(P)}$.


\section*{Agradecimentos}

Agradecemos apoio da FAPESP e do DM da UNESP - Câmpus de Bauru.

\end{document}